\documentclass[reqno,a4paper,12pt,russian,english]{amsart}
\usepackage{amsmath,amsthm,amsfonts,amssymb}
\usepackage{verbatim, graphicx, ifthen, enumitem}
\usepackage[T1]{fontenc}
\usepackage [applemac,koi8-r,latin9] {inputenc}
\usepackage{upquote,bbm}
\usepackage{babel}
\usepackage[bookmarksopen=false,pdfborder={0 0 0},pdfborderstyle={},colorlinks=true,linkcolor=black,citecolor=black,urlcolor=blue]{hyperref}

\allowdisplaybreaks

\let\oldmarginpar\marginpar
\renewcommand\marginpar[1]{\-\oldmarginpar[\raggedleft\footnotesize #1]%
{\raggedright\footnotesize #1}}
\theoremstyle{plain}

\newtheorem*{thm*}{Theorem}

\def\clap#1{\hbox to 0pt{\hss#1\hss}}

\usepackage{calrsfs}
\usepackage{graphicx}
\usepackage{color}
\usepackage{perpage}

\gdef\SetFigFontNFSS#1#2#3#4#5{} 

\theoremstyle{remark}
\newtheorem*{qst*}{Question}

\theoremstyle{definition}

\theoremstyle{remark}

\newtheorem*{remark*}{Remark}

\newcommand{\Z}{\mathbb{Z}}

\newcommand{\N}{\mathbb{N}}
\newcommand{\R}{\mathbb{R}}

\def\L{\Lambda}
\def\O{S}
\def\l{\lambda}

\def\N{\mathbb{N}}
\def\Z{\mathbb{Z}}

\def\R{\mathbb{R}}
\def\C{\mathbb{C}}
\def\F{\mathcal{F}}

\def\1{\mathbbm{1}}

\begin{document}

\title[Exponential frames for unbounded sets]{Exponential frames on unbounded sets}
\author{Shahaf Nitzan, Alexander Olevskii and Alexander Ulanovskii}

\maketitle

\begin{abstract} For every set $\O$ of finite measure in $\R$ we construct a discrete set
              of real frequencies $\L$ such that the exponential system $\{\exp(i\l t),\l\in\L\}$
              is a frame in $L^2(\O)$.
\end{abstract}

   \section{Introduction}    This note can be viewed as a continuation of our previous paper [NOU].
                   In [NOU] we constructed  "good" sampling sets for the Paley--Wiener spaces $PW_\O$ of entire $L^2(\R)-$functions  with bounded  spectrum $\O$ in $\R$.
                   This construction  is based on a result in [BSS] on existence of                    well-invertible  sub-matrices of large orthogonal matrices.   Recently, an important progress in the latter area has been made in [MSS].
              Based on this,  we prove existence of exponential frames in $L^2(\O)$, for every unbounded set $\O$  in $\R$ of finite measure.

                     Recall  that a system of vectors  $E=\{u_j\}$ is a frame in a Hilbert space $H$
                    if there are positive constants $a,A$ such that
$$
a\|h\|^2\leq \sum_{u_j\in E} |\langle h,u_j\rangle|^2\leq A\|h\|^2 \ \ \  \ \forall h\in H.
$$The numbers   $a$ and $A$ above are called  frame bounds.

                      Given a discrete set $\L$ in $\R$, we denote by $$
                                    E(\L):= \{e^{i\l t}\}_{\l\in \L}
$$ the system of exponentials with frequencies in $\L$.


                      Exponential frames $E(\L)$ in $L^2(\O)$  (equivalently, stable sampling sets $\L$ for $PW_\O$)  have been  carefully studied from
                    different points of view.    There is a large number of results in the area.              In the classical case when $\O$ is an interval, such systems were essentially   characterized  by  Beurling [B] in terms of the so-called "lower uniform density" of $\L$. A complete description of exponential frames for intervals is given by Ortega--Cerd\'{a} and  Seip [OS].
                     However, the problem of existence of  exponential  frames  for unbounded sets remained open.
                     The following result fills this gap by showing that for every set $\O$ of finite measure, the space $L^2(\O)$ admits an exponential frame:

                          \medskip\noindent{\bf Theorem 1} {\sl  There are positive constants $c,C$ such that for every set $\O\subset\R$ of finite measure there is a discrete set $\L\subset\R$ such that $E(\L)$  is a frame in $L^2(\O)$ with frame bounds $c |\O|$ and $C|\O|$.}

\medskip Here  by $|\O|$              we denote the measure of $\O$.

\medskip\noindent{\bf Remark 1}. The frame bounds are essential in many contexts, since they characterize the "quality"  of frame decompositions. Assume that an exponential system $E(\L)$ forms an orthogonal basis in $L^2(\O)$. One can easily check that in this case $E(\L)$ is a frame in $L^2(\O)$ with frame bounds $a=A=|\O|$. These are, in a sense, the "optimal" frame bounds. In general,  there may be no exponential orthogonal basis in $L^2(\O).$ However, Theorem~1 shows that an exponential frame in $L^2(\O)$ always exists with "almost" (up to  fixed multiplicative constants) optimal frame bounds.

              \medskip      \noindent{\bf Remark 2}.  A similar to Theorem 1 result regarding the existence of complete exponential systems  $E(\L)$ in $L^2(\O)$  (equivalently, existence of uniqueness sets $\L$ for $PW_\O$) is obtained in [OU] by an effective direct  construction.  That is not the case here, since the proof of Theorem A below in  [MSS] involves stochastic elements.

  \medskip      \noindent{\bf Remark 3}. Assume that $\O$ lies on an interval of length $2\pi d, d>0$. It follows from Lemma 10 below that a set $\L$ satisfying the conclusion of Theorem~1 can be chosen satisfying $\L\subset (1/d)\Z$.

 \medskip      \noindent{\bf Remark 4}. Assume that $\L$ satisfies the conclusions of Theorem~1. Then there are two absolute constants $k,K$ such that  the inequalities $$k|\O|< \frac{\#(\L\cap \Omega)}{|\Omega|}<K|\O|$$ hold whenever $\Omega$ is a sufficiently long interval in $\R$. In fact, one can choose any numbers  $k<1/2\pi$ and $K>4C$, where $C$ is the constant in Theorem 1. Then, as it was shown by Landau [L] (for a more elementary proof see [NO]), the left hand-side inequality above follows from the frame property of $E(\L)$.  The right hand-side inequality follows from Lemma 6 (ii) below.

                    \section{Well-invertible submatrices }

                     Our construction is based on the following result by Marcus, Spielman and Srivastava from [MSS]:

    \medskip\noindent{\bf Theorem A} {\sl Let $\epsilon>0$, and $u_1, . . . , u_m\in \C^n$ such that $\|u_i\|^2\leq\epsilon$ for all $i=1,...m$, and
\[
\sum_{i=1}^m| \langle w,u_i\rangle|^2=\| w\|^2 \ \ \ \ \ \forall w\in\C^n.
\]
Then there exists a partition of $\{1,...,m\}$ into $S_1$ and $S_2$, such that for each $j={1, 2}$,}
\begin{equation}\label{mss}
\sum_{i\in S_j}|\langle w,u_i\rangle|^2\leq\frac{(1+\sqrt{2\epsilon})^2}{2}\|w\|^2\qquad\forall w\in\C^n.
\end{equation}
\medskip

Observe that, clearly, $(1+\sqrt{2\epsilon})^2\leq 1+5\sqrt\epsilon$ when $\epsilon<1$.

\medskip\noindent{\bf Remark 5}. Let $\epsilon<1.$
 Since$$
\sum_{i\in S_1}|\langle w,u_i\rangle|^2 =\| w\|^2-\sum_{i\in S_2}|\langle w,u_i\rangle|^2, $$  estimate (1) shows that the two-sided estimate holds for each $j=1,2$:
\begin{equation}\label{1}
\frac{1-5\sqrt{\epsilon}}{2}\|w\|^2\leq \sum_{i\in S_j}|\langle w,u_i\rangle|^2\leq\frac{1+5\sqrt{\epsilon}}{2}\|w\|^2\ \ \ \qquad\forall w\in\C^n.
\end{equation}

   \medskip

The following corollary (see Corollary F.2 in [HO]) gives a reformulation of Theorem A in a form well prepared for an induction process:

                          \medskip\noindent{\bf Corollary B} {\sl
Let $v_1, . . . , v_k\in \C^n$
be such that $\|v_i\|^2\leq\delta$ for all $i=1,...,k$. If
\[
\alpha \|w\|^2\leq \sum_{i=1}^k|\langle w,v_i\rangle|^2\leq\beta \|w\|^2\qquad\forall w\in\C^n,
\]
with some numbers $\alpha>\delta$ and $\beta$, then there exists a partition of $\{1,...,k\}$ into $S_1$  and $S_2$ 
such that for each $j=1,2,$
 \begin{equation}\label{3}
\frac{1-5\sqrt{\delta/\alpha}}{2}\alpha \|w\|^2\leq \sum_{i\in S_j}|\langle w,v_i\rangle|^2\leq\frac{1+5\sqrt{\delta/\alpha}}{2}\beta \|w\|^2\qquad\forall w\in\C^n.
\end{equation}
}

\medskip\noindent  For the sake of completeness, we reproduce the proof.

Let $M:\C^n\rightarrow\C^n$ be the operator defined by $Mw=\sum_{i=1}^k\langle w, v_i\rangle v_i$. Observe that   $M$ is positive and that
$$
\alpha\|w\|^2\leq \| M^{1/2}w\|^2\leq \beta\|w\|^2 \ \ \ \ \ \forall w\in \C^n.
$$

Set $u_i=M^{-1/2}v_i.$ Then $\|u_i\|^2 \leq \|v_i\|^2/\alpha\leq\delta/\alpha.$  Further, for all $w\in\C^n$,
\[
\sum_{i=1}^k\langle w, u_i\rangle u_i=M^{-1/2}\sum_{i=1}^k\langle M^{-1/2}w, v_i\rangle v_i=M^{-1/2}MM^{-1/2}w=w.
\]
We see that $u_i$ satisfy the assumptions of Theorem A with $m=k$ and $\epsilon=\delta/\alpha<1$. Hence, there is a partition of $\{1,...,k\}$ into two sets $S_1$ and $S_2$ satisfying (\ref{1}). 
Using the right hand-side of (\ref{1}) we get
\[
\sum_{i\in S_j}|\langle w,v_i\rangle|^2=\sum_{i\in S_j}|\langle M^{1/2}w,u_i\rangle|^2 \leq\frac{1+5\sqrt{\epsilon}}{2}\|M^{1/2}w\|^2\leq
\]
\[
\leq\frac{1+5\sqrt{\epsilon}}{2}\beta\|w\|^2= \frac{1+5\sqrt{\delta/\alpha}}{2}\beta\|w\|^2.
\]
The proof of  the left hand-side of (\ref{3}) is similar.

\medskip

We will use an elementary lemma:

\medskip\noindent{\bf Lemma 1} {\sl Let $0<\delta<1/100$, and  let $\alpha_j, \beta_j, j=0,1,...,$ be defined inductively
$$
 \alpha_0= \beta_0=1,\  \alpha_{j+1}:=\alpha_j\frac{1-5\sqrt{\delta/\alpha_j}}{2},\ \ \beta_{j+1}:=\beta_j\frac{1+5\sqrt{\delta/\alpha_j}}{2}.
$$Then there exist a positive absolute constant $C$ and a number  $L\in\N$  such that}$$a_j\geq 100\delta, j\leq L, \ \ 25\delta\leq a_{L+1}<100\delta, \ \ b_{L+1}<Ca_{L+1}.$$

\medskip
\noindent{\bf Proof}. Clearly, if $a_j\geq 100\delta$ then $$
 \frac{\alpha_j}{4}\leq\alpha_{j+1}<\frac{\alpha_j}{2}.$$

Denote by  $L\geq 1$ the greatest  number such that $a_L\geq 100\delta,$ and set $\gamma_j:=5\sqrt{\delta/a_j}, j\leq L.$ Then $\gamma_{L-j}<2^{-1-j/2}$. It follows that
$$
\prod_{j=0}^L\frac{1+\gamma_j}{1-\gamma_j}<C:=\prod_{j=0}^\infty\frac{1+2^{-1-j/2}}{1-2^{-1-j/2}}.
$$This gives
$b_{L+1}<Ca_{L+1}$, and the lemma follows.

\medskip

We will need the following


\medskip\noindent{\bf Lemma 2} {\sl Assume the hypothesis of Theorem A are fulfilled and that  $\| u_i\|^2=n/m, i=1,...,m$. Then there is a subset $J\subset\{1,...,m\}$ such that \begin{equation}\label{2} c_0\frac{n}{m}\|w\|^2\leq \sum_{i\in J}|\langle w,u_i\rangle|^2\leq C_0\frac{n}{m}\|w\|^2\qquad\forall w\in\C^n, \end{equation}where $c_0$ and $C_0$ are some {\it absolute} positive constants.}

\medskip\noindent{\bf Proof}. 
 If $n/m\geq 1/100,$ then (\ref{2}) holds with $J=\{1,...,m\}$ and $C_0=c_0=100.$

Assume  $\delta:=n/m<1/100.$ Let $\alpha_j$ and $\beta_j$ be as defined in Lemma 1. Then the vectors $v_i=u_i$ satisfy the assumptions of Corollary B with $\alpha_0=\beta_0=1.$ Hence,  a set $J_1\subset\{1,...,m\}$ exists  such that
$$
\alpha_1\| w\|^2\leq \sum_{i\in J_1}|\langle w,u_i \rangle|^2\leq\beta_1\| w\|^2 \ \ \ \ \forall w\in \C^n.
$$

Since  $\alpha_1\geq \alpha_{L}>100\delta$, we may  apply Corollary~B the second time to get a set $J_2\subset J_1$ such that the two-sided inequality above holds with $J_2,\alpha_2$ and $\beta_2$, and so on. Since $\alpha_{L}>100\delta$,  Corollary~B can be applied  $L$ times.
         We thus obtain  a set $J_{L+1}\subset\{1,...,m\}$ for which the two-sided inequality holds with $\alpha_{L+1}$ and $\beta_{L+1}$.  From Lemma 1 it follows that (\ref{2})  is true with $J=J_{L+1}.$



\medskip
We now reformulate Lemma 1 in terms more convenient for our application.  
Given a matrix $A$ of order $m\times n$ and a subset $J\subseteq \{1,...,m\}$, we denote by $A(J)$ the
sub-matrix of $A$ whose rows belong to the index set $J$.

\medskip\noindent{\bf Lemma 3} {\sl
There exist positive constants $c_0, C_0>0$, such that whenever $A$ is an $m\times n$ matrix which is a sub-matrix of some
$m\times m$ orthonormal matrix, and such that all of its rows
have equal $l^2$ norm, one can find a subset $J\subset\{1,...,m\}$ such that
 \begin{equation}\label{mat-est}
c_0\frac{n}{m} \|w\|^2\leq \|A(J)w\|^2\leq C_0\frac{n}{m} \|w\|^2\qquad\forall w\in\C^n.
\end{equation}
}

\section{Auxiliary results}

In what follows we write $F=\hat f$, where  $f$ is the Fourier transform of $F$:
$$
f(x)=\frac{1}{\sqrt{2\pi}}\int_\R e^{-itx}F(t)\,dt.
$$

Given a discrete set $\L$, we denote by $d(\L)$ its separation constant
$$
d(\L):=\inf_{\l,\l'\in\L,\l\ne\l'}|\l-\l'|.
$$
Given a sequence of sets $\L_j$ satisfying $d(\L_j)\geq d>0$ for all $j$, a set $\L$ is called the weak limit
of $\L_j$ if for every $\epsilon>0$ and for every interval $\Omega = (a, b), a, b\not\in\L$,
both inclusions $\L_j\cap \Omega\subset (\L\cap \Omega)+(-\epsilon,\epsilon)$ and $\L\cap \Omega\subset (\L_j\cap \Omega)+(-\epsilon,\epsilon)$ hold for all but a finite number of
$j$'s. The standard diagonal procedure implies that if  $\L_j$ satisfy $d(\L_j)\geq d>0$ for all $j$, then there is a subsequence which weakly converges to some (maybe, empty) set $\L$ satisfying $d(\L)\geq d.$

Recall that the Paley--Wiener space $PW_\O$ is defined as the space of all functions $f\in L^2(\R)$ such that $\hat f$ vanishes a.e. outside $\O$. When the measure of $S$ is finite, we have$$\int_{S}|F(t)|\,dt\leq \|F\|\sqrt{|S|} \ \  \ \ \ \ \ \ \forall F\in L^2(S).$$Here $\|F\| $ means the $L^2-$norm of $F$. Hence,  $\hat f\in L^1(\R)$ for every $f\in PW_S$, and so every function $f\in PW_S$ is continuous.

Sometimes it will be more convenient for us to work with the Paley--Winer space $PW_\O,$ rather than $L^2(\O).$ In this connection we observe that by taking the Fourier transform, Theorem 1 is equivalent to the following statement:

{\sl There exist positive constants $c,C$ such that for every set $\O\subset\R, |\O|<\infty,$
there is a discrete set $\L\subset\R$ such that }
\begin{equation}\label{0}
c|\O|\|f\|^2\leq \sum_{\l\in\L} |f(\l)|^2\leq C|\O|\|f\|^2 \:\:\:\:\:\:\:\:\:\:\:\:\forall
f\in PW_\O
\end{equation}

\medskip
We will prove (\ref{0}) with the constants $C=C_0$ and $c=c_0/(36C_0)$, where $c_0$ and $C_0$ are the constants in Lemma 3.

\medskip
We will  need the following Bessel's inequality (see [Y], Ch. 4.3): Given a set $\L$ satisfying $d(\L)>0$ and a bounded set $\O$,  there is a constant $K$ which depends only on $d(\L)$ and the diameter of $\O$  such that
$$
\sum_{\l\in\L}|f(\l)|^2\leq K\|f\|^2 \ \ \ \ \ \ \ \ \  \forall f\in PW_\O.
$$

The proof of Theorem 1 below uses three auxiliary lemmas:

\medskip\noindent{\bf Lemma 4} {\sl Let $\O$ be a bounded set of positive measure and let $\L_k\subset\R$ be a sequence of  sets satisfying $d(\L_k)>\delta>0, k=1,2,...$, which converges weakly to some set $\L$. Then$$
\lim_{k\to\infty}\sum_{\l\in\L_k} |f(\l)|^2=\sum_{\l\in\L} |f(\l)|^2\:\:\:\:\:\:\:\:\:\:\:\:\forall
f\in PW_{\O}.
$$}

\noindent{\bf Proof}. Take any function $f\in PW_{\O}$, and pick up a point  $x_l\in  [l\delta-\delta/2,l\delta+\delta/2]$ such that
$$|f(x_l)|=\max_{|x-l\delta|\leq\delta/2} |f(x)|  \ \ \  \ \ \ \ \forall l\in\Z.$$ Since $x_{l+2}-x_l\geq\delta$,  the sequence $x_k$ is a union of two sets each having separation constant $\geq\delta.$ By Bessel's inequality,  we see that
$$\sum_{k\in\Z}|f(x_k)|^2<\infty.
$$

Let $R>0$, and write
$$
\big|\sum_{\l\in\L_k} |f(\l)|^2-\sum_{\l\in\L} |f(\l)|^2\big|\leq \big|\sum_{\l\in\L_k, |\l|<R} |f(\l)|^2-\sum_{\l\in\L, |\l|<R} |f(\l)|^2\big|+ 2\sum_{|k|\geq R/\delta}|f(x_k)|^2.
$$The first term in the right hand-side tends to zero as $k\to\infty$ whenever $\pm R\not\in\L$, while the second one tends to zero as $R\to\infty.$ This proves the lemma.

\medskip\noindent{\bf Lemma 5} {\sl
Let $\O_1\subseteq \O_2\subseteq...$ be an increasing sequence of bounded sets in $\R$ with $\O=\cup_k\O_k$ being a set of finite measure. Let $\L\subset\R,d(\L)>0,$ and positive $k,K$ be such that  the inequalities
\begin{equation}\label{j}
k\|f_j\|^2\leq \sum_{\l\in\L} |f_j(\l)|^2\leq K\|f_j\|^2 \:\:\:\:\:\:\:\:\:\:\:\:\forall
f_j\in PW_{\O_j}
\end{equation}hold for every $j$.  Then
\begin{equation}\label{jj}k\|f\|^2\leq \sum_{\l\in\L} |f(\l)|^2\leq K\|f\|^2 \:\:\:\:\:\:\:\:\:\:\:\:\forall
f\in PW_{\O}.\end{equation}}

\medskip\noindent{\bf Proof}.
Given a function $f\in PW_\O$, let $f_j\in PW_{\O_j}$ be the Fourier transform of the function $\hat f\cdot 1_{S_j}$, where $1_{S_j}$ is the indicator function of $S_j$. Then the $L^1-$norm of $\hat f-\hat{f_j}$ tends to zero as $j\to\infty,$ and so the functions $f_j(x)$ converge uniformly to $f(x)$.

  For every $R>0$ we have,
\[
\sum_{\l\in\L, |\l|<R}|f_j(\l)|^2\leq K\|f_j\|^2.
\]
Taking the limit as $j\to\infty$, we obtain
\[
\sum_{\l\in\L, |\l|<R}|f(\l)|^2\leq  K\|f\|^2.
\]
By letting $R\to\infty$, we obtain the right hand-side inequality in (\ref{jj}). Using this inequality, we get
\[
\big(\sum_{\l\in\L}|f(\l)|^2\big)^{1/2}\geq \big(\sum_{\l\in\L}|f_j(\l)|^2\big)^{1/2}- \big(\sum_{\l\in\L}|(f-f_j)(\l)|^2\big)^{1/2}\geq
\]
\[
k^{1/2}\|f_j\|- K^{1/2}\|f-f_j\|^2.
\]
Taking the limit as $j\to\infty$, we prove the left hand-side inequality in (\ref{jj}).

\medskip\noindent{\bf Lemma 6} {\sl Assume that the inequality
\begin{equation}\label{h}
\sum_{\l\in\L}|f(\l)|^2\leq C |\O|\|f\|^2 \ \ \ \ \ \forall f\in PW_\O
\end{equation}is true for some $C>0, \O\subset\R, |\O|<\infty,$ and   $\L\subset\R$. Then

(i) There is a constant $\eta>0$ which depends only on $\O$ such that $$\#(\L\cap \Omega)\leq 9C,$$for every interval $\Omega\subset\R, |\Omega|=\eta$.

(ii) There is a constant $K>0$ which depends only on $\O$ such that$$\frac{\#(\L\cap \Omega)}{|\Omega|}\leq 4C|\O|,$$for every interval $\Omega\subset\R, |\Omega|\geq K$. }

\medskip\noindent{\bf Proof}. (i)
Denote by $h\in PW_\O$ the Fourier transform of the indicator function $1_\O$. Then $h(x)$ is continuous,
$$
h(0)=\frac{|\O|}{\sqrt{2\pi}}, \ \|h\|^2= \|1_\O\|^2=|\O|.
$$Choose $\eta>0$ so small that $|h(x)|> |\O|/3, |x|\leq \eta/2$.  Then, applying (\ref{h}) for $f=h$, we see that the statement (i) of Lemma 6 holds for $\Omega=[-\eta/2,\eta/2].$ To complete the proof, it suffices to observe that every function $h(x-x_0), x_0\in\R,$ belongs to $PW_\O$.

(ii) Take any function $g\in PW_\O$ satisfying $\|g\|=1$, and choose a number $R$  such that
$$
\int_{-R}^R|g(x)|^2\,dx\geq \frac{1}{2}.$$ Assume $K>2R$. We now   apply (\ref{h}) to the function $f(x):=g(x-s)$ and integrate   over $(-K,K)$ with respect to $s$:$$\int_{-K}^K\sum_{\l\in\L}|g(\l-s)|^2\,ds\leq 2KC |\O|.
$$
When $|\l|<K/2$, we have $$
\int_{-K}^K|g(\l-s)|^2\,ds\geq \int_{-R}^R|g(s)|^2\,ds\geq\frac{1}{2}.
$$
We conclude that
$$
\frac{\#(\L\cap(-K/2,K/2))}{2}\leq\int_{-K}^{K}\sum_{\l\in\L}|g(\l-s)|^2\,ds\leq 2KC |\O|.
$$This proves statement (ii).

\section{Proof of Theorem 1}


The proof of Theorem 1 will consist of a series of lemmas.

\medskip\noindent
{\bf Lemma 7} {\sl
Let $n,m\in\N, n<m$. For every set
$$
\O=\bigcup_{r\in I}\left[\frac{2\pi r}{m},\frac{2\pi(r+1)}{m}\right],\ I\subset\{0,...,m-1\}, \# I=n,
$$
there is a set $\L\subset\Z$ such that \begin{equation}\label{00}
c_0|\O|\|f\|^2\leq \sum_{\l\in\L} |f(\l)|^2\leq C_0|\O|\|f\|^2 \:\:\:\:\:\:\:\:\:\:\:\:\forall
f\in PW_\O,
\end{equation} where $c_0,C_0$ are the constants in Lemma 3.}

\medskip
\noindent
{\bf Proof}.
Observe that $|\O|=2\pi n/m,$ and denote by $$\F_I:=(e^{i\frac{2\pi jr}{m}})_{r\in
I,j=0,...,m-1 }$$ the submatrix of the Fourier matrix $\F$ whose
columns are indexed by $I$. Since the matrix $(\sqrt m)^{-1}\F$ is
orthonormal, by Lemma 3 there exists $J\subset\{0,...,m-1\}$ such that
\begin{equation}\label{f}
c_0n\Vert w\Vert^2\leq\Vert \F_I(J)w\Vert^2_{l_2(J)}\leq C_0n\Vert w\Vert^2,
\:\:\:\:\:\:\:\:\:\:w\in l_2(I).\end{equation}

 Observe that every function $F\in L^2(\O)$ can be
written as $$F(t)=\sum_{r\in I}F_r(t-\frac{2\pi r}{m}), $$ where
$F_r\in L^2(0,\frac{2\pi}{m})$ is defined by
$$F_r(t):=F(t+\frac{2\pi r}{m}){\bf 1}_{[0,\frac{2\pi}{m}]}(t).$$ Therefore, every function $f\in PW_\O$ admits a representation
$$
f(x)=\sum_{r\in I}e^{i\frac{2\pi
r}{m}x}f_r(x),\:\:\:\:\:\:\:\:\:\: f_r\in PW_{[0,
\frac{2\pi}{m}]},
$$
where the functions $e^{i\frac{2\pi r}{m}x}f_r(x)$ are orthogonal
in $L^2(\R)$. We note that for every function $h\in PW_{[0, 2\pi/m]}$
we have,
\begin{equation}\label{ogb}
\frac{2\pi}{m}\|h\|^2=\sum_{\l\in m\Z} |h(\l)|^2.
\end{equation}

 We now verify that the sequence
$$
\L:=\{j+km: j\in J, k\in\Z\}
$$
satisfies (\ref{00}). Take any function $f\in PW_\O$. Then
$$
\sum_{j\in J}\sum_{k\in\Z}|f(j+km)|^2=\sum_{j\in
J}\sum_{k\in\Z}\left|\sum_{r\in I}e^{i\frac{2\pi rj
}{m}}f_r(j+km)\right|^2.
$$
  For every $j\in J$ we apply (\ref{ogb}) to the function $\sum_{r\in I}e^{i\frac{2\pi rj
}{m}}f_r(x)$. We find that the last expression is equal to
$$
\frac{2\pi }{m}\sum_{j\in J}\int_{\R}\left|\sum_{r\in
I}e^{i\frac{2\pi rj }{m}}f_r(x)\right|^2dx= \frac{2\pi
}{m}\int_{\R}\|\F_I(J)(f_r(x))_{r\in I}\|^2_{l_2(J)}dx.$$ By
inequality (\ref{f}) we have on one hand,
$$
\sum_{\l\in\L}|f(\l)|^2\geq c_0\frac{n }{m}\int_{\R}\sum_{r\in
I}|f_r(x)|^2dx=
$$$$c_0\frac{n }{m}\int_{\R}\sum_{r\in I}|e^{i\frac{2\pi
r}{m}x}f_r(x)|^2dx= c_0\frac{n }{m}\int_{\R}|\sum_{r\in
I}e^{i\frac{2\pi r}{m}x}f_r(x)|^2dx=$$$$c_0\frac{n
}{m}\int_{\R}|f(x)|^2dx,
$$
while on the other hand, applying the same computation, we get
$$
\sum_{\l\in\L}|f(\l)|^2\leq C_0\frac{n }{m}\int_{\R}\sum_{r\in
I}|f_r(x)|^2dx=
 C_0\frac{n
}{m}\int_{\R}|f(x)|^2dx.
$$

This completes the proof.

\medskip\noindent{\bf Lemma 8} {\sl  For every compact set $\O\subset[0,2\pi]$ of positive measure there is a set $\L\subset\Z$ such that (\ref{00}) holds.}

\medskip
This follows immediately from Lemma 7, since every such set $\O$ can be covered by a set from Lemma 7 whose measure is arbitrarily close to $|\O|$.

\medskip\noindent{\bf Lemma 9} {\sl  For every set $\O\subset[0,2\pi]$ of positive measure there is a set  $\L\subset\Z$ such that (\ref{00}) holds.}

\medskip\noindent{\bf Proof}. It suffices to prove Lemma 9 for open sets $\O$.  Let $\O$ be such a set end
let $\O_1\subset \O_2\subset...$ be an increasing sequence of compact sets such that $\O = \cup_j \O_j$. By Lemma 8, there exist sets $\L_j\subset\Z$ such that \begin{equation}\label{000}
c_0|\O|\|f_j\|^2\leq \sum_{\l\in\L_j} |f_j(\l)|^2\leq C_0|\O|\|f_j\|^2 \:\:\:\:\:\:\:\:\:\:\:\:\forall
f_j\in PW_{\O_j},
\end{equation}where $c_0,C_0$ are the constants in Lemma 3. Since $PW_{\O_j}\subset PW_{\O_k}, k> j,$ we have
\begin{equation}\label{l}
c_0|\O_k|\|f_j\|^2\leq \sum_{\lambda\in \L_k}|f_j(\l)|^2\leq C_0|\O_k|\|f_j\|^2 \ \ \ \ \ \forall f_j\in PW_{\O_j}
\end{equation}
We may assume that $\L_k$ converge weakly to some set $\L\subset\Z$.  Using Lemma 4, we take the limit as $k\to\infty$:
\begin{equation}\label{01}
c_0|\O|\|f_j\|^2\leq \sum_{\lambda\in \L}|f_j(\l)|^2\leq C_0|\O|\|f_j\|^2 \ \ \ \ \ \forall f_j\in PW_{\O_j}.
\end{equation}
Now, the result follows from Lemma 5.

\medskip\noindent{\bf Lemma 10} {\sl  For every bounded set $\O$ of positive measure there is a set  $\L\subset (1/d)\Z$ such that (\ref{00}) holds, where $d$ is any positive number such that $\O$ lies on an interval of length $2\pi d$.}

\medskip Observe that the translations of $\O$  change neither the frame property of $E(\L)$ nor the frame constants. So, it suffices to assume that $\O\subset[0,2\pi d]$. Then the result follows from Lemma 9 by re-scaling.

\medskip\noindent{\bf Proof of Theorem 1}. We may assume that $\O$ is an unbounded set of finite measure.

Let  $\O_1\subset \O_2\subset...$ be any sequence of bounded sets satisfying $\O=\cup_j \O_j$.   By Lemma 10, there exist discrete sets $\L_j$ such that (\ref{000}) is true.
Since $PW_{S_j}\subset PW_{S_k},j<k$, we see that  (\ref{l}) holds for all $j<k$.

 By Lemma 6 (i),  there is a number $\eta>0$ and an integer $r$ which depends only on the constant $C_0$ in (\ref{0}) (it is easy to check that one may take $r\leq 36 C_0$) such that every set $\L_k$ can be can be splitted up into $r$ subsets $\L_k^{(l)}$ satisfying
$d(\L_k^{(l)})\geq \eta, l=1,...,r.$ By taking an appropriate subsequence, we may assume that each $\L_k^{(l)}$ converges weakly to some set $\L^{(l)}$ as $k\to\infty.$
By Lemma 4,  we may take  limit in (\ref{l}) as $k\to\infty$:
\[
c_0 |\O|\|f_j\|^2\leq \sum_{k=1}^r\sum_{\l\in\L^{(k)}}|f_j(\l)|^2\leq C_0|\O|\|f_j\|^2\qquad\forall f_j\in PW_{\O_{j}}.
\]

Set $\L:=\cup_{k=1}^r\L^{(k)}$. It may happen that the sets $\L^{(k)}$ have common points. Anyway, we have
$$
\sum_{\l\in\L}|f_j(\l)|^2\leq\sum_{k=1}^r\sum_{\l\in\L^{(k)}}|f_j(\l)|^2\leq r\sum_{\l\in\L}|f_j(\l)|^2.
$$   From the latter inequalities, it readily  follows that
\[
\frac{c_0}{r} |\O|\|f_j\|^2\leq \sum_{\l\in\L}|f_j(\l)|^2\leq C_0|\O|\|f_j\|^2\qquad\forall f_j\in PW_{\O_{j}}.
\]
Theorem 1 now follows easily from Lemma 5.

\section{References}

\noindent [BSS] J. Batson, D. A. Spielman, and N. Srivastava. Twice--Ramanujan sparsifiers. SIAM Rev. 56 (2014), no. 2, 315--334.

\medskip
\noindent [B] A. Beurling. Balayage of Fourier--Stieltjes Transforms. In: The collected
Works of Arne Beurling. Vol. 2, Harmonic Analysis. Birkhauser, Boston 1989.

\medskip
\noindent [HO]	N. J. A. Harvey, N. Olver. Pipage rounding, pessimistic estimators and matrix concentration. Proc. of the Twenty-Fifth Annual ACM-SIAM Symposium on Discrete Algorithms (2014),  926--945. ISBN: 978-1-611973-38-9.

\medskip\noindent [L] H. J. Landau. Necessary density conditions for sampling and interpolation of certain entire functions. Acta
Mathematica (1967), v. 117, 37--52.

\medskip\noindent [MSS]
A. Marcus, D. A. Spielman, and N. Srivastava. Interlacing families II: Mixed characteristic
polynomials and the Kadison-Singer problem, June 2013. arXiv:1306.3969.

\medskip\noindent [NOU] S. Nitzan, A. Olevskii, A. Ulanovskii.  A few remarks on sampling of signals with small spectrum. Proceeding of the Steklov Institute of Mathematics (2013), v. 280, 240--247.

\medskip\noindent [NO] S. Nitzan, A. Olevskii,
Revisiting Landau's density theorems for Paley--Wiener spaces. C. R. Mathematique (2012),
v. 350, no. 9--10, 509--512.

\medskip\noindent [OU] A. Olevskii, A. Ulanovskii. Uniqueness sets for unbounded spectra. C. R. Math. Acad. Sci. Paris 349 (2011), no. 11--12, 679--681.

\medskip\noindent [OS]
J. Ortega--Cerd\'{a}, K. Seip. Fourier frames. Annals of Mathematics (2002), v. 155 (3), 789--806.

\medskip\noindent
[Y] R.M: Young. An introduction to Nonharmonic Fourier Series. Academic Press. 2001.

\end{document}